\documentclass[a4paper]{article}
\pdfoutput=1

\usepackage{IEK10_Single} \usepackage{natbib}
\usepackage{orcidlink}

\def\IEK10{
  Forschungszentrum Jülich GmbH,
  Institute of Energy and Climate Research,
  Energy Systems Engineering (IEK-10),
  Jülich 52425,
  Germany
}
\def\RWTH{
  RWTH Aachen University
  Aachen 52062,
  Germany
}
\def\JARA{
  JARA-ENERGY,
  Jülich 52425,
  Germany
}
\def\SVT{
  RWTH Aachen University,
  Process Systems Engineering (AVT.SVT),
  Aachen 52074,
  Germany
}
\def\TUB{
    Technical University of Berlin, 
    Berlin 10623,
    Germany
}

\newcommand{\mytitle}{
Normalizing Flow-based Day-Ahead Wind Power Scenario Generation for Profitable and Reliable Delivery Commitments by Wind Farm Operators
}

\newcommand{\affil}{
  \begin{itemize}[leftmargin=3mm, itemsep=0mm]
    \item[$^a$]\IEK10
    \item[$^b$]\RWTH
    \item[$^c$]\TUB
    \item[$^d$]\JARA
    \item[$^e$]\SVT
\end{itemize}
}

\def\firstAuthor{Eike Cramer}
\newcommand{\myauthor}{
\firstAuthor$^{a,b}$\orcidlink{0000-0002-6882-5469}, 
Leonard Paeleke$^{a,c}$,
Alexander Mitsos$^{d,a,e}$\orcidlink{0000-0003-0335-6566}, 
Manuel Dahmen$^{a,*}$\orcidlink{0000-0003-2757-5253} }
\newcommand{\correspondingauthor}{Manuel Dahmen, \IEK10 \newline E-mail: m.dahmen@fz-juelich.de}
\author{\myauthor}

\begin{document}

  \thispagestyle{firststyle}

  \begin{center}
    \begin{large}
      \textbf{\mytitle}
    \end{large} \\
    \myauthor
  \end{center}

  \vspace{0.5cm}

  \begin{footnotesize}
    \affil
  \end{footnotesize}

  \vspace{0.5cm}

  \begin{abstract}
    We present a specialized scenario generation method that utilizes forecast information to generate scenarios for day-ahead scheduling problems. 
In particular, we use normalizing flows to generate wind power scenarios by sampling from a conditional distribution that uses wind speed forecasts to tailor the scenarios to a specific day. 
We apply the generated scenarios in a stochastic day-ahead bidding problem of a wind electricity producer and analyze whether the scenarios yield profitable decisions. 
Compared to Gaussian copulas and Wasserstein-generative adversarial networks, the normalizing flow successfully narrows the range of scenarios around the daily trends while maintaining a diverse variety of possible realizations. 
In the stochastic day-ahead bidding problem, the conditional scenarios from all methods lead to significantly more stable profitable results compared to an unconditional selection of historical scenarios. 
The normalizing flow consistently obtains the highest profits, even for small sets scenarios.   \end{abstract}
    
  \vspace{0.5cm}
  
  \noindent 
  \textbf{Keywords}: Wind power; scenario generation; stochastic programming; stability 
  \section{Introduction}
Since the liberalization of the electricity markets, electricity is traded on the electricity SPOT markets \citep{EPEX2021Documentation,mayer2018electricity}.
The auction-based format of the day-ahead markets requires electricity producers and large-scale consumers to specify fixed amounts of electricity they want to buy or sell one day prior to the delivery \citep{EPEX2021Documentation}. 
Thus, renewable electricity producers have to account for the uncertain and non-dispatchable nature of renewable electricity generation from wind and photovoltaic when submitting their bids \citep{perez2012impacts,mayer2018electricity,mitsos2018challenges}.

To find profitable solutions, operators often leverage optimization techniques from the process systems engineering (PSE) community \citep{ZHANG2016114,grossmann_2021advanced}.
In particular, scheduling optimization identifies cost-optimal operational setpoints and leverages variable electricity prices \citep{schafer2020wavelet,leo2021stochastic}.
To address the uncertainty stemming from the uncertain renewable electricity production and volatile price curves, scheduling problems are often implemented as stochastic programs that include the uncertainty in the problem formulation \citep{conejo2010decision,grossmann_2021advanced}.
Typically, stochastic programs are based on scenarios, e.g., possible realizations of renewable production trajectories \citep{conejo2010decision,morales2013integrating,chen2018advances}. 
The PSE community has been at the forefront of finding solutions to scheduling problems and stochastic programs for decades \citep{grossmann1978optimum,halemane1983optimal,pistikopoulos1995novel,sahinidis2004optimization}. 
Thus, energy system scheduling problems are solved successfully by the PSE community \citep{ZHANG2016114,schafer2019model,schafer2020wavelet}.
Many PSE examples address electricity procurement for power-intensive processes and demand-side-management \citep{ZHANG2016114,zhang2016risk, leo2021stochastic}.
Examples with energy focus include \cite{garcia2008stochastic}, who derive a stochastic bidding problem for a wind producer with pumped hydro storage, and \cite{liu2015bidding}, who propose a model to obtain bidding curves for a micro-grid considering distributed generation. In their book, \cite{conejo2010decision} derive a wind producer bidding problem considering both price and production uncertainties.

While most works focus on optimization problem formulations and their solutions, obtaining high-quality scenarios is also critical for operational success. 
The scenarios for stochastic programming either stem from historical data or specialized scenario generation methods \citep{conejo2010decision}. 
Established methods for scenario generation often utilize univariate, i.e., step-by-step prediction, approaches like classical autoregressive models \citep{sharma2013wind} or autoregressive neural networks \citep{vagropoulos2016ann,voss2018Residential}.
As opposed to univariate models, multivariate modeling techniques model a series of time steps in a single prediction step. This makes them particularly suitable for day-ahead operation problems as the multivariate predictions best capture the correlations throughout the day \citep{ziel2018day} and can be set up to model the distribution of the given time horizon, i.e., the time frame between 00:00\,am and 11:59\,pm of the following day. 
Prominent multivariate scenario generation models are Gaussian copulas \citep{pinson2009probabilistic,staid2017generating,camal2019scenario} as well as deep generative models like generative adversarial networks (GANs) \citep{chen2018model,jiang2018scenario,wei2019short} and variational autoencoders (VAEs) \citep{zhanga2018optimized}. 

Despite their widespread application, the training success of both GANs and VAEs is sometimes poor and their loss functions are difficult to judge as they are not directly concerned with the quality of the generated data \citep{salimans2016improved, borji2018pros}. Furthermore, GANs and VAEs often result in a mode collapse, i.e., the models converge to a single feasible scenario instead of describing the true probability distribution \citep{arjovsky2017principled}.
Besides GANs and VAEs, normalizing flows are another type of deep generative model \citep{papamakarios2021normalizing}. The major advantage of normalizing flows is their training via direct log-likelihood maximization, which leads to interpretable loss functions and stable convergence \citep{rossi2018mathematical}.
In prior works, normalizing flows performed well for multivariate probabilistic time series modeling \citep{rasul2021multivariate}, for scenario generation of residential loads \citep{zhang2019scenario,ge2020modeling} and wind and photovoltaic electricity generation \citep{dumas2021deep,cramer2022pricipalcomponentflow}.

Many authors argue that their scenario generation approach samples high-quality scenarios \citep{pinson2009probabilistic, chen2018model, zhang2019scenario}. However, a connection of scenario generation to downstream applications in stochastic programming is missing in most contributions. 
Exceptions are \cite{zhanga2018optimized} and \cite{wei2019short} who both solve operational problems for wind-solar-hydro hybrid systems.
However, their respective VAE and GAN are restricted to unconditional scenario generation, i.e., they sample unspecific scenarios without considering the day-ahead setting and without including forecasts or other available information. 
For a day-ahead bidding problem, this can potentially lead to suboptimal solutions based on an unrealistic scenario set containing many unlikely scenarios.
Meanwhile, conditional scenario generation incorporates forecasts and other available information to specifically tailor the scenarios to the following day.
Examples of conditional scenario generation are the Gaussian copula approach by \cite{pinson2009probabilistic} and the normalizing flow by \cite{dumas2021probabilistic}, where only \cite{dumas2021probabilistic} solve a stochastic optimization problem using quantiles derived from the conditional normalizing flow. 

\cite{Kaut2003Evaluation} derive two criteria to evaluate scenario generation methods for stochastic programming. 
First, \cite{Kaut2003Evaluation} define \textit{stability} via the variance in the optimal objective values, where different sets of scenarios sampled from the same scenario generation method should yield similar optimal objectives in the stochastic program.
Second, \cite{Kaut2003Evaluation} define the \textit{bias} of a scenario generation method as the difference between the optimal objective of the scenario-based stochastic program and the optimal objective obtained with the true distribution. 
While stability can be evaluated by solving multiple instances of the same stochastic program based on different scenario sets, the bias of a scenario generation method is impossible to evaluate for the day-ahead scenario generation problem as the true distribution is unknown and cannot be approximated via historical scenarios. 
Notably, none of the previously published works on DGM-based scenario generation test for the criteria defined by \cite{Kaut2003Evaluation}.

This work extends our previous work on normalizing flow-based scenario generation \citep{cramer2022pricipalcomponentflow} to perform conditional scenario generation \citep{zhang2019scenario,dumas2021probabilistic} of wind power generation with wind speed forecasts as conditional inputs, i.e., we use the wind speed forecast to generate day-ahead wind power generation scenarios that are specifically tailored to the given day.
We then apply the generated scenarios in a stochastic day-ahead wind electricity producer bidding problem based on \cite{garcia2008stochastic} and \cite{conejo2010decision}.
We compare the results obtained using the normalizing flow scenarios with unconditional historical scenarios and two other multivariate conditional scenario generation approaches, namely, the well-established Gaussian copula \citep{pinson2009probabilistic} and the recently very popular Wasserstein-GAN (W-GAN) \citep{chen2018model}.
Our analysis shows that all conditional scenario generation methods result in significantly more profitable decisions compared to the historical data and that the profits obtained using the normalizing flow scenarios are the highest among all considered methods. 
Unlike \cite{wei2019short} or \cite{dumas2021probabilistic}, we also perform a statistical investigation of the stability defined by \cite{Kaut2003Evaluation}. 
In particular, we consider that most stochastic programs can only be solved for small sets of scenarios, which makes stability increasingly difficult to achieve. Hence, we solve the stochastic problem for limited scenario-set-sizes to investigate their applicability to stochastic programs that cannot be solved for a high number of scenarios.

The remainder of this work is organized as follows:
Section~\ref{sec:Cond_NormFlow} details the concept of conditional density modeling using normalizing flows.
Then, Section~\ref{sec:Cond_dayaheadScheme} details the conditional day-ahead scenario generation method and reviews the input-output relation of normalizing flows, Gaussian copulas, and W-GANs.
Section~\ref{sec:Cond_CaseStudyData} draws a comparison of historical scenarios and scenarios generated using the three different methods based on the analysis outlined in \cite{pinson2012evaluating} and \cite{cramer2021validation}.
Section~\ref{Sec:Cond_CaseStudy2} introduces the stochastic bidding problem and analyzes the stability and the obtained profits of the different scenario sets. 
Finally, Section~\ref{sec:Cond_Conclusion} concludes this work.   \section{Conditional density estimation using normalizing flows}\label{sec:Cond_NormFlow}
Normalizing flows are data-driven, multivariate probability distribution models that use invertible neural networks $T: \mathbb{R}^D \rightarrow \mathbb{R}^D$ to describe a data probability density function (PDF) as a change of variables of a $D$-dimensional Gaussian distribution \citep{kobyzev2019normalizing, papamakarios2021normalizing}:
\begin{equation*}
\begin{aligned}
    \mathbf{x} &= T(\mathbf{z})\\
    \mathbf{z} &= T^{-1}(\mathbf{x})
\end{aligned}
\end{equation*}
Here, $\mathbf{x}\in X\subset \mathbb{R}^D$ are samples of the data and $\mathbf{z}\in \mathcal{N}(\mathbf{0}_D, \mathbf{I}_D)$ are the corresponding vectors in the Gaussian distribution with $\mathbf{I}_D$ being the $D$-dimensional identity matrix. 
New data $\mathbf{x}$ is generated by drawing samples $\mathbf{z}$ from the known Gaussian distribution and transforming them via the forward transformation $T(\cdot)$. 
Since the transformation between the data and the Gaussian is a change of variables, the PDF of the data can be expressed explicitly via the inverse transformation $T^{-1}(\cdot)$ using the change of variables formula \citep{papamakarios2021normalizing}:
\begin{equation}
    p_X(\mathbf{x}) = \phi(T^{-1}(\mathbf{x}))  \left| \det \mathbf{J}_{T^{-1}}(\mathbf{x}) \right| \label{Eq:ChangeOfVariablesInverse}
\end{equation}
Here, $\mathbf{J}_{T^{-1}}$ is the Jacobian of the inverse transformation $T^{-1}$, and $p_X$ and $\phi$ are the PDFs of the data and the Gaussian, respectively. 
Intuitively, Equation~\eqref{Eq:ChangeOfVariablesInverse} describes a projection of the data onto the Gaussian and a scaling of the distribution's volume to account for the constant probability mass.
If the transformation $T$ is a trainable function, the normalizing flow can be trained via direct log-likelihood maximization using the log of Equation~\eqref{Eq:ChangeOfVariablesInverse}.

To describe a conditional PDF $p_{X\vert Y}(\mathbf{x}\vert \mathbf{y})$ with conditional inputs $\mathbf{y}\in Y$, i.e., the joint PDF of $X$ and $Y$ where the realization $\mathbf{y}\in Y$ is known, the transformation $T$ and its inverse $T^{-1}$ must accept the conditional information vector $\mathbf{y}$ in addition to the transformed variables $\mathbf{z}$ and $\mathbf{x}$, respectively \citep{winkler2019learning}:
\begin{equation*}
\begin{aligned}
    \mathbf{x} &= T(\mathbf{z}, \mathbf{y})\\
    \mathbf{z} &= T^{-1}(\mathbf{x}, \mathbf{y})
\end{aligned}
\end{equation*}
If $T$ remains differentiable for any fixed value of the conditional inputs $\mathbf{y}$, the likelihood can still be described using the change of variables formula:
\begin{equation}
    p_{X\vert Y}(\mathbf{x}\vert \mathbf{y}) = \phi(T^{-1}(\mathbf{x},\mathbf{y}))  \left| \det \mathbf{J}_{T^{-1}}(\mathbf{x},\mathbf{y}) \right| \label{Eq:ChangeOfVariablesConditionalInverse}
\end{equation}

\begin{figure*}
    \centering
    \includegraphics[width=\textwidth]{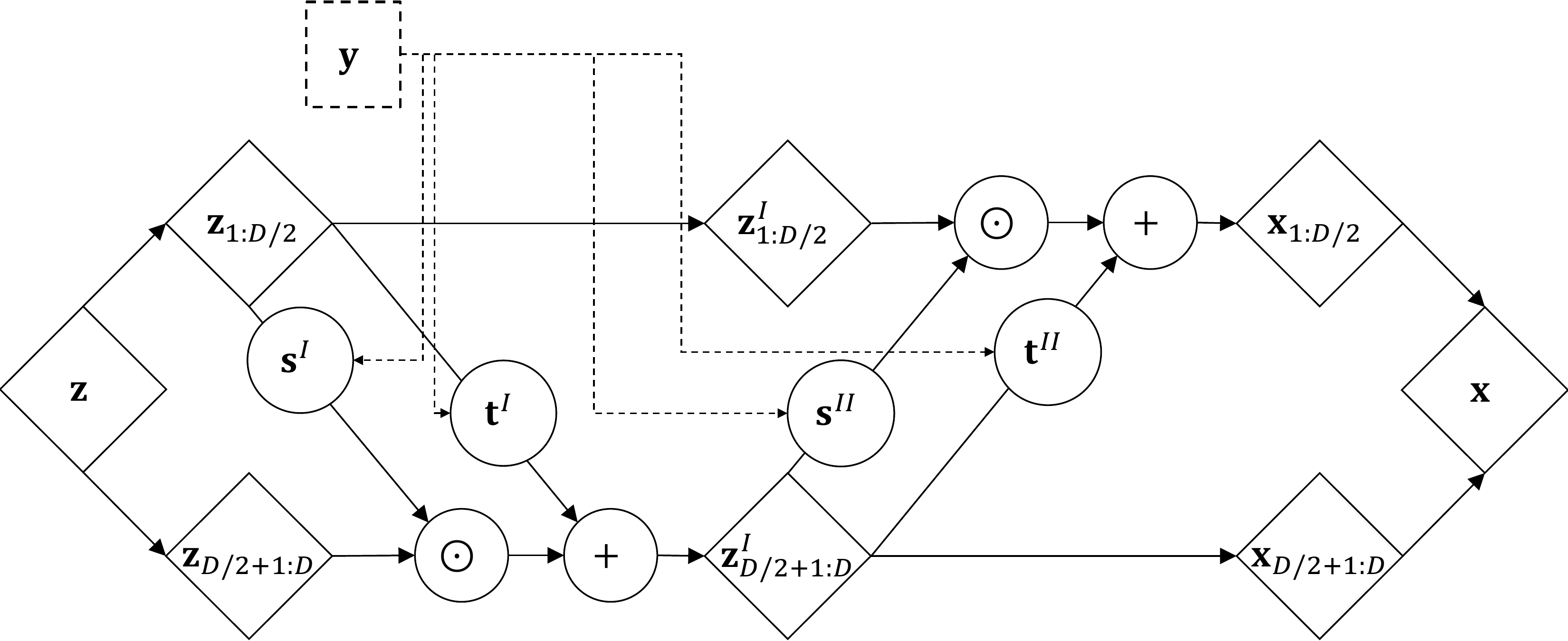}
    \caption{Example RealNVP architecture containing two conditional affine coupling layers \citep{dinh2016realNVP,winkler2019learning}, with
    conditioner models $\mathbf{s}^{I}$, $\mathbf{t}^{I}$, $\mathbf{s}^{II}$, and $\mathbf{t}^{II}$, 
    Gaussian sample vector $\mathbf{z}$,
    data sample vector $\mathbf{x}$,
    intermediate sample vector $\mathbf{z}^{I}$,
    conditional input vector $\mathbf{y}$.
    The indices $1:D/2$ and $D/2+1:D$ refer to the two halves of the data vectors, respectively. The dashed lines indicate the flow of the conditional input data $\mathbf{y}$.
    }
    \label{fig:CondRealNVP}
\end{figure*}
In this work, we employ the real non-volume preserving transformation (RealNVP) \citep{dinh2016realNVP}, which is based on a composition of affine coupling layers. In each coupling layer, one half of the data vector undergoes an affine transformation, which is parameterized via functions of the remaining half of the data vector:
\begin{equation}
    \begin{aligned}
        \mathbf{x}_{1:D/2}   =& \mathbf{z}_{1:D/2}  \\
        \mathbf{x}_{D/2+1:D} =& \mathbf{z}_{D/2+1:D}\odot \exp(\mathbf{s}(\mathbf{z}_{1:D/2},\mathbf{y})) 
         + \mathbf{t}(\mathbf{z}_{1:D/2},\mathbf{y})
    \end{aligned}
    \label{Eq:RealNVPConditional}
\end{equation}
Here, $\odot$ denotes element-wise multiplication, 
the indices $1:D/2$ and $D/2+1:D$ refer to the two halves of the data vectors, respectively, 
and $\mathbf{s}$ and $\mathbf{t}$ are the so-called conditioner models. 
Notably, the partial transformation of the data vectors keeps the inputs of $\mathbf{s}$ and $\mathbf{t}$ identical for both the forward and the inverse transformation. Thus, the conditioner models can be implemented as any standard feed forward neural network \citep{dinh2016realNVP}. 
Furthermore, the clever design of the affine coupling layer results in lower triangular Jacobians. Hence, the Jacobian determinant required for the likelihood computation is simply given by the product over the diagonal elements. The log-form Jacobian determinant used for training then is:
\begin{equation*}
    \log \det \mathbf{J}_{\text{RealNVP}}(\mathbf{z}) = \sum_{k = D/2+1}^{D} \mathbf{s}(\mathbf{z}_{1:D/2},\mathbf{y})
\end{equation*}
Large and highly expressive normalizing flows can be built using compositions of Equation~\eqref{Eq:RealNVPConditional} in an alternating manner. Figure~\ref{fig:CondRealNVP} shows an illustrative sketch of an exemplary RealNVP architecture with two conditional affine coupling layers.

In \cite{cramer2022pricipalcomponentflow}, we showed that normalizing flows sample uncharacteristically noisy scenarios when applied to sample for the distributions of renewable electricity time series, due to their inherent lower-dimensional manifold structure.
To address the issue, we proposed dimensionality reduction based on the principal component analysis (PCA). 
In this work, we use PCA to reduce the dimensionality of the data $\mathbf{x}$ and the Gaussian samples $\mathbf{z}$. The conditional input vectors $\mathbf{y}$ are not affected by the PCA.
For more information on the effects of manifolds we refer to \cite{brehmer2020flows}, \cite{behrmann2021understanding}, and \cite{cramer2022pricipalcomponentflow}.   \section{Day-ahead scenario generation}\label{sec:Cond_dayaheadScheme}
This work addresses scenario generation with a particular focus on applications in day-ahead scheduling problems. Thus, all scenarios describe a possible realizations covering the time between 00:00\,am and 11:59\,pm of the following day.
In particular, we generate day-ahead wind power generation scenarios and use day-ahead forecasts of wind speeds as conditional inputs to narrow down the range of possible trajectories and make the scenarios specific to the following day. For reference, we include a comparison to historical data, which represents unconditional scenarios, i.e., randomly drawn sampled from the full distribution $p_X(\mathbf{x})$ that does not consider the wind speed forecasts.
Meanwhile, the scenario generation methods aim to fit models of the full conditional PDF $p_{X\vert Y}(\mathbf{x}\vert \mathbf{y})$ that are valid for every possible wind power realization $\mathbf{x} \in X$ and every possible day-ahead wind speed forecast $\mathbf{y} \in Y$. 
In the application, the wind speed predictions are known one day prior to the scheduling horizon and the scenario generation models are evaluated for fixed conditional inputs $\mathbf{y}$.
Normalizing flows, Gaussian copulas, and W-GANs all employ multivariate modeling approaches, i.e., the models generate full daily trajectories in a vector form \citep{pinson2009probabilistic, ziel2018day, chen2018model}. 
All models use multivariate Gaussian samples $\mathbf{z}$ and the wind speed forecast vectors, i.e., the conditional information $\mathbf{y}$, as inputs to generate wind power generation scenario vectors $\mathbf{x}$. 
For a given fixed wind speed forecast $\mathbf{y} = \text{const.}$, sampling and transforming multiple Gaussian samples $\mathbf{z}$ results in a set of wind power generation scenarios, i.e., the Gaussian acts as a source of randomness to generate sets of scenarios instead of point forecasts: 
\begin{equation*}
    \mathbf{x}_i = T(\mathbf{z}_i,\mathbf{y}=\text{const.})\quad \forall i \in 1, \dots, \# \text{Scenarios}
\end{equation*}
Here, $T(\cdot)$ can be any of the scenario generation models. 
For more details on the evaluation of the different models, we refer to our supplementary material and the papers by \cite{pinson2009probabilistic} and \cite{chen2018model}.

All models generate capacity factor scenarios, i.e., the actual production scaled to installed capacity, of the 50 Hertz transmission grid in the years 2016 to 2020 \citep{DataSource}. The year 2019 is set aside as a test set to avoid complications in the stochastic programming case study due to the unusual prices resulting from the COVID-19 pandemic \citep{narajewski2020changes,badesa2021ancillary}. 
To avoid including test data in the scenario sets, the unconditional historical scenarios are drawn from the training set.
The 15\,min recording interval renders 96-dimensional scenario vectors that fit the 24\,h time horizon of a day-ahead bidding problem. 
The day-ahead wind speed forecasts have hourly resolution and are obtained from the reanalysis data set ``Land Surface Forcings V5.12.4'' of MERRA-2 \citep{globalmodelingassimilationoffice2015} which is based on previously recorded historical data.
We use the predictions at the coordinates 53.0$^\circ$\,N, 13.0$^\circ$\,E, in the center of the 50 Hertz region.
Note that due to potential wind speed forecast errors and agglomeration effects in the power generation, there is no direct known functional relationship between the wind speed forecast and the realization of regionally distributed power generation.

Due to numerically singular Jacobians and non-invertible transformations \citep{behrmann2021understanding,cramer2022pricipalcomponentflow}, full-space normalizing flows fail to accurately describe the distribution of daily wind time series trajectories residing on lower-dimensional manifolds \citep{cramer2022pricipalcomponentflow}.
Therefore, we use PCA \citep{pearson1901pca} to reduce the data dimensionality following our recent contribution \citep{cramer2022pricipalcomponentflow}. 
We select the number of principal components based on the explained variance ratio, i.e., the amount of information maintained by the PCA \citep{pearson1901pca}. For an explained variance ratio of 99.95\%, we obtain 18 principal components to represent the original 96-dimensional scenario vectors. 
The adversarial training algorithm for the W-GAN \citep{arjovsky2017wasserstein} did not converge consistently for the considered learning problem. Thus, the results presented below are drawn from the best performing model out of 20 different trained models w.r.t. the metrics outlined in Section~\ref{sec:Cond_CaseStudyData}.
For more detailed information on the implementation, we refer to the supplementary material.    \section{Conditional wind power scenario generation}\label{sec:Cond_CaseStudyData}
We start by analyzing the scenarios without a specific application in mind. To this end, we present some examples, analyze the described distributions, and investigate whether the models can identify the correct daily trends. 

\begin{figure*}
    \centering
    \includegraphics[width=\textwidth]{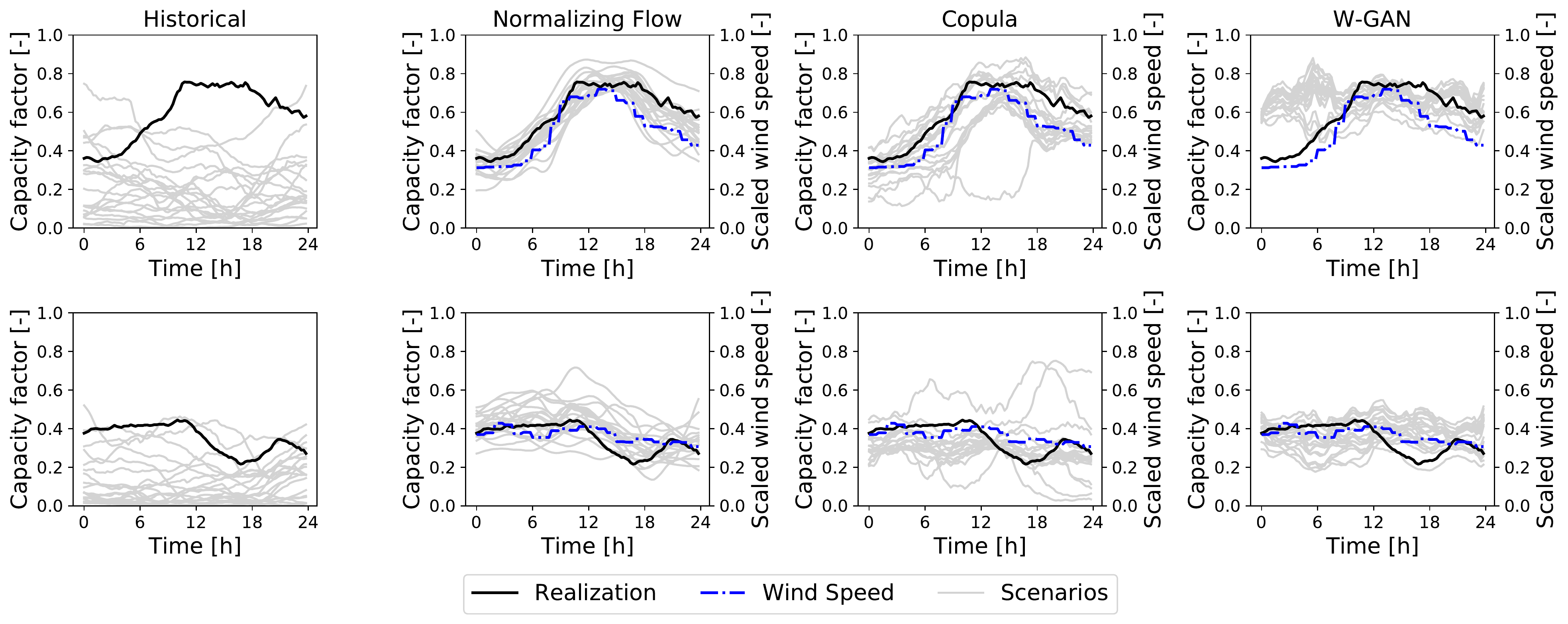}
    \caption{
         20 wind capacity factor scenarios (``Scenarios'') each from the historical scenario set (left), normalizing flow (center-left), Gaussian copula (center-right), and W-GAN (right) in relation to the realized wind capacity factor (``Realization''). 
         The plots for the the scenario generation methods include the conditional input (``Wind Speed''). Realization and scenarios on the left y-axis and scaled wind speed predictions on right y-axis.
        Top: March 29th, 2019, bottom: September 15th, 2019. Data from 50 Hertz region \citep{DataSource}.}
    \label{fig:Cond_ExampleScenario}
\end{figure*}
Figure~\ref{fig:Cond_ExampleScenario} shows example scenarios for two randomly selected days of the test year 2019.
The left, center-left, center-right, and right columns show historical scenarios and scenarios sampled from the conditional normalizing flow, the Gaussian copula, and the W-GAN, respectively.
The historical scenarios are randomly selected from the training set and are, therefore, unspecific to the respective days. Thus, they fail to identify the daily trends and show large discrepancies for both days. 
For both example days in Figure~\ref{fig:Cond_ExampleScenario} the normalizing flow identifies and follows the general trend of the realized wind capacity factor. For the presented examples, the realization lies within the span of the scenarios. 
Similarly, the Gaussian copula also identifies the trend of the realization. 
However, there are some scenarios with significantly higher or lower capacity factors in the case of both days, i.e., the Gaussian copula appears prone to sample outliers that do not follow the trend. 
The W-GAN-generated scenarios fail to identify the trend and, instead, appear tightly agglomerated and only represent the daily average of the realization, which can be observed for the morning hours of the first day and, to a lesser extend, on the afternoon hours of the second day.
The failed identification of the trend is likely due to a mode collapse of the W-GAN, which is a frequently observed phenomenon with GANs \citep{arjovsky2017principled}. Mode collapse happens when the adversarial training algorithm converges to a small range of realistic scenarios but fails to identify the true distribution. Note that due to the multivariate modeling approach of generating full daily trajectories, this type of deviation may occur at any time step throughout the day.

\begin{figure*}
    \centering
    \includegraphics[width=\textwidth]{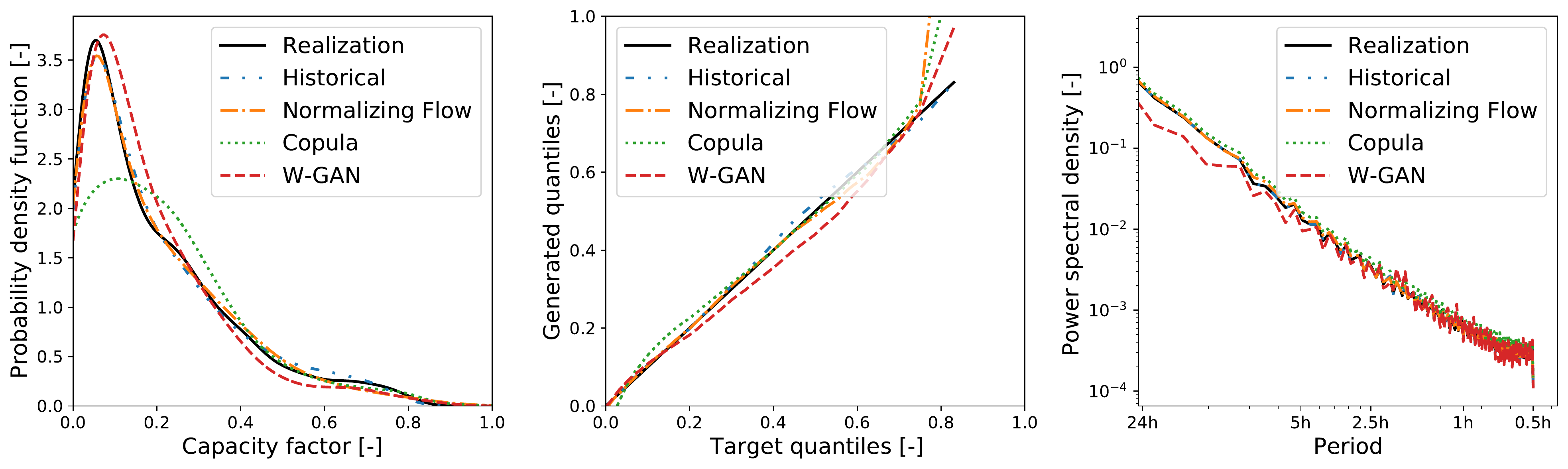}
    \caption{
        Distribution and fluctuational behavior of generated wind capacity factor scenarios from historical training data (``Historical''), normalizing flow (``Normalizing Flow''), Gaussian copula (``Copula''), and W-GAN (``W-GAN'') in comparison to historical test data (``Realization'') \citep{DataSource}.
        Left: marginal probability density function (PDF) estimated using kernel density estimation \citep{parzen1962estimation},
        center: quantile-quantile plots (Q-Q-plots) \citep{chambers2018graphical},
        right: power spectral density (PSD) estimated using Welch's method \citep{welch1967use}.
    }
    \label{fig:Cond_PDF_QQ_PSD}
\end{figure*}

To gain insight into the quality of the full scenario sets, we analyze whether the scenario generation methods are able to reproduce the probability distributions and the frequency behavior of the actual time series by looking at the full year of 2019 in comparison to the eventual realization. 
To this end, we look into the marginal PDF \citep{parzen1962estimation}, the quantile distribution in Q-Q plots \citep{chambers2018graphical}, and the power spectral density (PSD) \citep{welch1967use}. For a detailed introduction to the interpretation of PDF and PSD, we refer to our previous work \citep{cramer2021validation}.
Figure~\ref{fig:Cond_PDF_QQ_PSD} shows the marginal PDFs (left), the Q-Q plots (center), and the PSD (right) of the historical- and the generated scenarios from the normalizing flow, the Gaussian copula, and the W-GAN in comparison to the realizations in 2019. 

In Figure~\ref{fig:Cond_PDF_QQ_PSD}, the historical scenarios and the normalizing flow scenarios describe the test set PDF well and show good matches of the quantile distribution in the Q-Q-plot.
Meanwhile, the Gaussian copula produces a broader PDF with a much lower peak than the realization, whereas the W-GAN's PDF shows a shift towards higher values. 
The Q-Q-plot also shows the shift of the W-GAN generated distribution, as there is an offset between the W-GANs and the other quantile lines. 
The poor distribution match by the Gaussian copula is likely due to the linear quantile regression that is unable to represent the nonlinear relation between the predicted wind speed and the capacity factor. Furthermore, the copula relies on linear interpolation of quantiles which can inflate the PDF in the tails and, thus, lead to outlier sampling \citep{pinson2009probabilistic}.
The W-GAN can theoretically model any distribution \citep{goodfellow2014generative}. In our analysis, however, the adversarial training algorithm was very difficult to handle with the time series data and often resulted in poor fits.
The presented results are the best of 20 training runs in terms of matching the criteria in Figure~\ref{fig:Cond_PDF_QQ_PSD}. 
Meanwhile, both the Gaussian copula and the normalizing flow with PCA converge consistently and typically yield the presented results after the first training attempt.

The Q-Q-plot reveals that all methods yield distributions with longer tails than the realizations, i.e., they produce scenarios with higher capacity factors than the maximum realized capacity factor. 
The reason is that for days with the highest capacity factor of the year, even higher capacity factors are still feasible as the realizations never reach the full installed capacity. Also, the log tail of the PDF makes the offset appear inflated as it only occurs for the 99-th and 100-th percentile. 
Note that both Copula and W-GAN are restricted to sample from the [0,1] interval via the boundaries of the inverse CDF \citep{pinson2009probabilistic} and the tanh output activation function, respectively. 
Meanwhile, the normalizing flow has no such restriction and yields some scenarios surpassing 1, which leads to the normalizing flow having the strongest deviation in the Q-Q-plot. Although these scenarios are theoretically infeasible, they have a very low probability and can efficiently be removed in postprocessing. 

The PSD in Figure~\ref{fig:Cond_PDF_QQ_PSD} shows a good match of the frequency behavior by the historical scenarios, the normalizing flow, and the Gaussian copula. 
The W-GAN is close to the overall power-law, i.e., the slope of the PDF curve, of the data, but fails to match the exact frequency behavior.

In addition to the analysis of the full scenario sets in Figure~\ref{fig:Cond_PDF_QQ_PSD}, we also compute the energy score (ES) for each day in 2019. 
ES is a quantitative measure for the assessment of multivariate scenario generation models that compares the conditional scenario set with the respective realization \citep{gneiting2008assessing,pinson2012evaluating}:
\begin{equation*}
    \text{ES} = 
        \frac{1}{N_S} \sum_{s=1}^{N_S} \vert\vert \mathbf{x} - \hat{\mathbf{x}}_s \vert\vert_2
        - \frac{1}{2{N_S}^2} \sum_{s=1}^{N_S} \sum_{s'=1}^{N_S} \vert\vert \hat{\mathbf{x}}_s - \hat{\mathbf{x}}_{s'} \vert\vert_2
\end{equation*}
Here, $\mathbf{x}$ is the realization vector, $\hat{\mathbf{x}}_s$ are the scenario vectors, $N_S$ is the number of scenarios, and $\vert\vert \cdot \vert\vert_2$ is the 2-norm.
The energy score is a negatively oriented score, i.e., lower values indicate better results.
The two parts of the energy score reward closeness to the realization and diversity of the scenario set, respectively.

\begin{figure*}
    \centering
    \includegraphics[width=\textwidth]{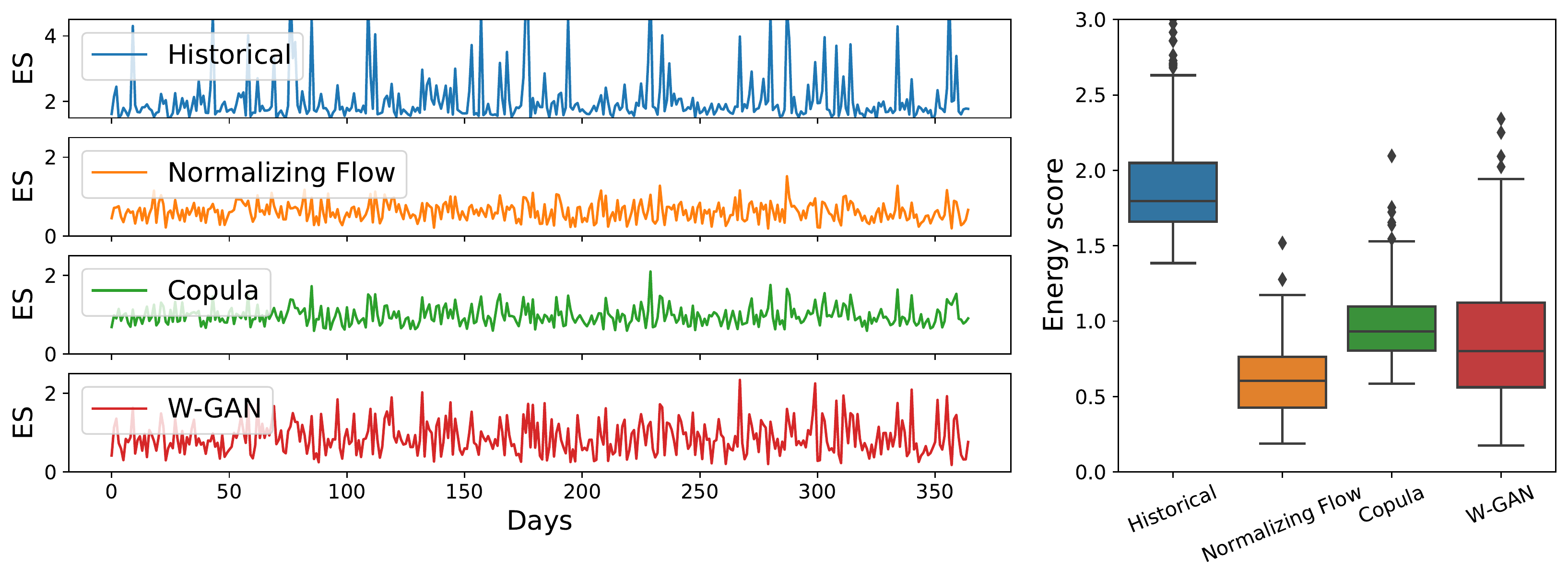}
    \caption{
    Energy score (ES) \citep{gneiting2008assessing, pinson2012evaluating} over all days in 2019 (left) and boxplots (right).
    Historical- (``Historical'') and generated scenarios from normalizing flow (``Normalizing Flow''), Gaussian copula (``Copula''), and Wasserstein-GAN (``W-GAN'').
    Boxes indicate quartiles and diamonds indicate outliers \citep{waskom2021seaborn}.
    Note different y-scale for historical ES. 
    }
    \label{fig:Cond_EnergyScore}
\end{figure*}

In Figure~\ref{fig:Cond_EnergyScore}, we display the energy score for each day in 2019 as well as boxplots that showcase the overall energy score distributions for the historical data and the three different models. 
The normalizing flow energy score is lower on average compared to the Gaussian copula and the W-GAN, indicating a better fitting of the realizations and more diverse scenarios. The Gaussian copula shows the highest energy score which is likely a result of the outliers observed in Figure~\ref{fig:Cond_ExampleScenario}.
Furthermore, the normalizing flow leads to a narrow distribution of energy scores with few outliers, indicating consistently good results. 
Meanwhile, the historical scenarios consistently result in significantly higher energy scores compared to the conditional day-ahead scenario generation methods. 
By design, the unconditional historical scenarios do not identify the daily trends and are not generated specifically for the respective days. Thus, the deviations from the realizations penalized by the energy score are significant for most days.

In conclusion, we find that the conditional normalizing flow presented in Section~\ref{sec:Cond_NormFlow} generates scenarios that match the true distribution of realizations closely, while also providing a diverse set of possible realizations. 
Furthermore, the normalizing flow outperforms the Gaussian copula and the W-GAN with respect to all important metrics.
The historical scenarios describe the overall distribution well, but are not specific to the individual days and, hence, return poor results in day-ahead problem-specific metrics like the energy score.

   \section{Day-ahead bidding strategy optimization}\label{Sec:Cond_CaseStudy2}
We apply the scenarios generated in the previous section in a wind producer bidding problem based on \cite{garcia2008stochastic} and \cite{conejo2010decision}. We first state the problem formulation and then analyze the stability of the scenario generation methods for different numbers of scenarios based on the criterion defined by \cite{Kaut2003Evaluation}.
Finally, we investigate the obtained profits based on the different scenario sets.

\subsection{Wind producer problem formulation}
\begin{figure}
    \centering
    \includegraphics[width=0.7\columnwidth]{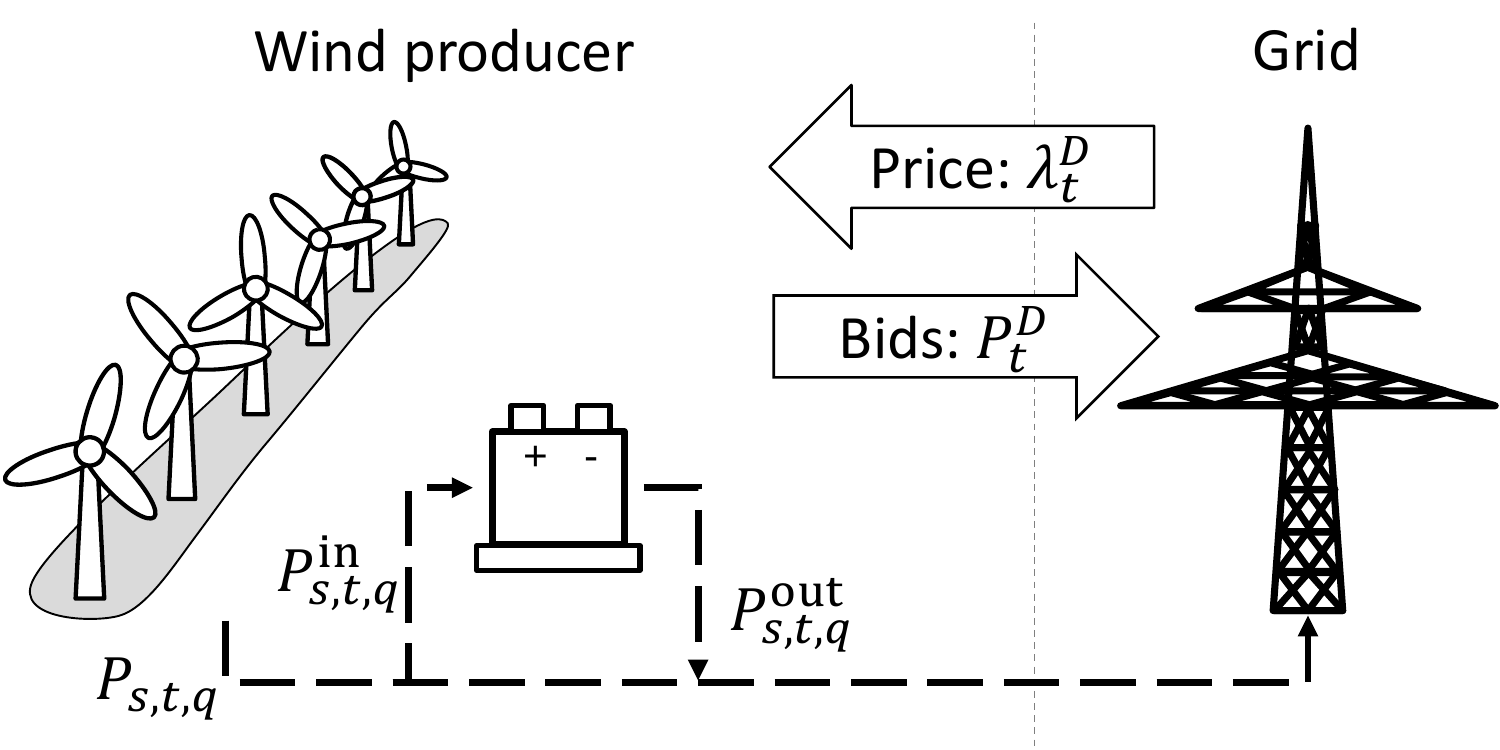}
    \caption{
        Structural setup of the wind producer problem from \cite{garcia2008stochastic} and \cite{conejo2010decision} with generated electricity $P_{s,t,q}$, (dis-) charging rates $P_{s,t,q}^{\text{in}}$ and $P_{s,t,q}^{\text{out}}$, placed bids $P_t^{D}$, and day-ahead electricity prices $\lambda_t^D$. The indices $s$, $t$, and $q$ indicate scenarios, hourly time intervals, and quaterhourly time intervals, respectively.
    }
    \label{fig:Cond_WindProducer}
\end{figure}

We consider the deterministic equivalent formulation \citep{birge2011introduction} of the stochastic wind producer problem from \cite{garcia2008stochastic} and \cite{conejo2010decision} shown in Figure~\ref{fig:Cond_WindProducer} that aims to find an optimal bidding schedule for the operator of a wind farm participating in the European Power Exchange (EPEX SPOT) market \citep{EPEX2021Documentation}. 
First, the operator places bids at the day-ahead auction market one day prior to delivery and, thereby, commits to deliver a certain amount of electricity $P^{D}_t$ during the given trading time interval $t$. 
The revenue made is given by $\lambda^{D}_t P^{D}_t \delta_h$, where $\lambda^{D}_t$ is the day-ahead price and $\delta_h =1\,h$ is the trading interval. 
As wind electricity generation is stochastic and non-dispatchable \citep{conejo2010decision}, the placed bids may not always be met by the actual production. 
To balance the difference between the placed bids and the actual production we allow for a small electricity storage that can store the electricity of up to 15\,min of maximum production. 
For any remaining production imbalance, we enforce a penalty on the absolute value of the imbalance \citep{garcia2008stochastic}. 
The full objective then reads \citep{garcia2008stochastic}:
\begin{equation}
     \underset{P^{D}_t}{\max} \sum_{t=1}^{N_T} \left[\lambda^{D}_t P^{D}_t \delta_h 
        - \omega \vert \lambda^{D}_{t} \vert \sum_{s=1}^{N_S}\pi_s \vert \Delta_{s,t} \vert \right] \label{Eq:Cond_WindProducerObjective}
\end{equation}
Here, $\Delta_{t,s}$ is the imbalance at time point $t$ and scenario $s$. The penalty term is based on the absolute values of the day-ahead price $\vert \lambda^{D}_{t} \vert$ to compensate for possible negative electricity prices \citep{garcia2008stochastic}.

The penalty term in Equation~\eqref{Eq:Cond_WindProducerObjective} contains the absolute value operator $\vert \cdot \vert$, leading to a nonlinear problem. 
However, any positive deviation can be avoided via curtailment of the plant and the imbalance will only take negative values in practice, which makes the absolute value operator obsolete.
Thus, the absolute imbalance is substituted by its negative parts to obtain a linear problem \citep{conejo2010decision}.
The complete linear formulation of the wind producer market participation problem including the electricity storage is shown in Problem~\eqref{Prob:Cond_WindProducerProblem}.

\begin{equation}
    \tag{WP}
    \label{Prob:Cond_WindProducerProblem}
    \begin{aligned}
\underset{P^{D}_t, P^{\text{in}}_{s,t,q}, P^{\text{out}}_{s,t,q}}{\max} \quad &\sum_{t=1}^{N_T} \left[\lambda^{D}_t P^{D}_t \delta_h 
        - \omega \vert \lambda^{D}_{t} \vert \sum_{s=1}^{N_S}\pi_s \Delta^{-}_{s,t} \right] \\
\text{s.t.}\quad 
& \text{SOC}_{s,t,q} =  \text{SOC}_{s,t,q-1} 
          + \eta \delta_q P^{\text{in}}_{s,t,q} - \frac{1}{\eta} \delta_q P^{\text{out}}_{s,t,q}, 
           &\forall s \in \mathcal{S}, \forall t\in \mathcal{T}, \forall q\in \mathcal{Q} \\
& \text{SOC}_{s, t=24, q=4} = \text{SOC}_{0}, &\forall s \in \mathcal{S}\\
& 
\Delta^{-}_{s,t}  \leq 
             \delta_h P^{D}_{t} - \delta_q \sum_{q\in \mathcal{Q}} P_{s,t,q} - \left(P^{\text{out}}_{s,t,q} - P^{\text{in}}_{s,t,q} \right),
              &\forall s \in \mathcal{S}, \forall t\in \mathcal{T} \\
& 0\leq P^{D}_{t} \leq P^{D,\max}, &\forall t\in \mathcal{T}\\
        & 0\leq \Delta^{-}_{s,t}, & \forall s \in \mathcal{S}, \forall t\in \mathcal{T} \\ 
        & 0\leq P^{\text{in}}_{s,t,q} \leq P^{\max}, &\forall s \in \mathcal{S}, \forall t\in \mathcal{T}, \forall q\in \mathcal{Q} \\
        & 0\leq P^{\text{out}}_{s,t,q} \leq P^{\max}, &\forall s \in \mathcal{S}, \forall t\in \mathcal{T}, \forall q\in \mathcal{Q} \\
        & 0\leq \text{SOC}_{s,t,q} \leq \text{SOC}^{\max}, &\forall s \in \mathcal{S}, \forall t\in \mathcal{T}, \forall q\in \mathcal{Q} \\
& \mathcal{S} = \{1,\dots, N_{S} \} \\
        & \mathcal{T} = \{1,\dots, N_T\} \\
        & \mathcal{Q} = \{1,\dots, 4 \} 
    \end{aligned}
\end{equation}
Tables~\ref{tab:Cond_Indices}, \ref{tab:Cond_Parameters}, and \ref{tab:Cond_Variables} list the indices, parameters, and variables of Problem~\eqref{Prob:Cond_WindProducerProblem}, respectively. 
Problem~\eqref{Prob:Cond_WindProducerProblem} is the deterministic equivalent of a two-stage stochastic program \citep{birge2011introduction}, where the delivery commitments are the first stage decisions and, the second stage decisions are the actual delivery and the storage operation.

The problem is implemented in pyomo \citep{pyomo}, version 6.2, and solved using gurobi \citep{gurobi}, version 9.5.
\begin{table}
    \centering
    \caption{Indices in Problem~\eqref{Prob:Cond_WindProducerProblem}. }
    \begin{tabularx}{\columnwidth}{lX}
        \hline
        Indices & Description           \\ 
        \hline
        $q$     & Quater hour interval  \\
        $s$     & Scenarios             \\ 
        $t$     & Hour interval         \\
        \hline
    \end{tabularx}
    \label{tab:Cond_Indices}
\end{table}

\begin{table}
    \centering
    \caption{Parameters in Problem~\eqref{Prob:Cond_WindProducerProblem}.}
    \begin{tabularx}{\columnwidth}{lXX}
        \hline
        Parameter           & Description                       & Value/Unit\\ 
        \hline 
        $\delta_h$          & Trading time interval             & 1\,h      \\
        $\delta_q$          & Production time interval          & 15\,min   \\
        $\eta$              & (Dis-) Charging efficiency        & 0.91      \\
        $\lambda^{D}_t$     & Day-Ahead Price                   & [EUR/MWh] \\
        $\omega$            & Penalty factor                    & 1.5       \\
        $\pi_s$             & Probability of scenario s         & $1/N_S$   \\
        $N_T$               & Number of time steps              & 24        \\
        $N_S$               & Number of scenarios               & [-]       \\
        $P^{D,\max}$        & Maximum production capacity       & 100\,MW   \\
        $P^{\max}$          & Maximum (dis-) charging rate      & 12.5\,MW  \\
        $P_{s,t,q}$         & Actual production                 & [MW]      \\
        $\text{SOC}^{\max}$ & Maximum battery capacity          & 25\,MWh   \\
        $\text{SOC}_{0}$    & Initial battery state of charge   & 12.5\,MWh \\
        \hline
    \end{tabularx}
    \label{tab:Cond_Parameters}
\end{table}

\begin{table}
    \centering
    \caption{Variables in Problem~\eqref{Prob:Cond_WindProducerProblem}.}
    \begin{tabularx}{\columnwidth}{lXl}
        \hline
        Variable             & Description & Unit                   \\
        \hline 
        $P^{D}_{t}$          & Bid at day-ahead market      & [MW]  \\
        $P^{in}_{s,t,q}$     & Charging rate                & [MW]  \\
        $P^{out}_{s,t,q}$    & Discharging rate             & [MW]  \\
        $\text{SOC}_{s,t,q}$ & Battery state of charge      & [MWh] \\
        $\Delta^{-}_{s,t}$   & Negative production imbalance& [MW] \\
        \hline
    \end{tabularx}
    \label{tab:Cond_Variables}
\end{table} Note that in Problem~\eqref{Prob:Cond_WindProducerProblem}, simultaneous charging and discharging of the storage is feasible, however, does not occur at the optimum due to the losses associated with using the storage.
The problem operates on both the trading time scale with hourly intervals and the production time scales with 15\,min intervals.

\subsection{Stability}
\cite{Kaut2003Evaluation} define stability and bias as criteria for the quality of a scenario generation method for stochastic programming. 
A scenario generation method is considered stable when different instances of the stochastic program based on different generated scenario sets result in similar optimal objective values. 
A small bias is achieved if the optimal objective value of the scenario-based formulation is close to the optimal objective value obtained by solving the stochastic program with the true distribution of the uncertain parameter. 
For the case of day-ahead scenario generation methods, the biases of the scenario generation methods are impossible to assess, as the true distribution for the wind power generation of the individual days is unknown and cannot be represented via historical data.

Problem~\ref{Prob:Cond_WindProducerProblem} is linear and can be solved efficiently. However, larger mixed-integer problems, as well as non-convex stochastic programs, often cannot be solved for a large number of scenarios due to the large computational effort \citep{birge2011introduction} and, small scenario sets must suffice for the stochastic program. Consequently, high stability becomes increasingly important for small numbers of scenarios. 
In the following, we solve 50 instances of Problem~\eqref{Prob:Cond_WindProducerProblem} for each day of 2019 and each scenario generation method using small scenario sets of only 3, 5, 10, 20, and 50 scenarios, each.

\begin{figure}
    \centering
    \includegraphics[width=0.7\columnwidth]{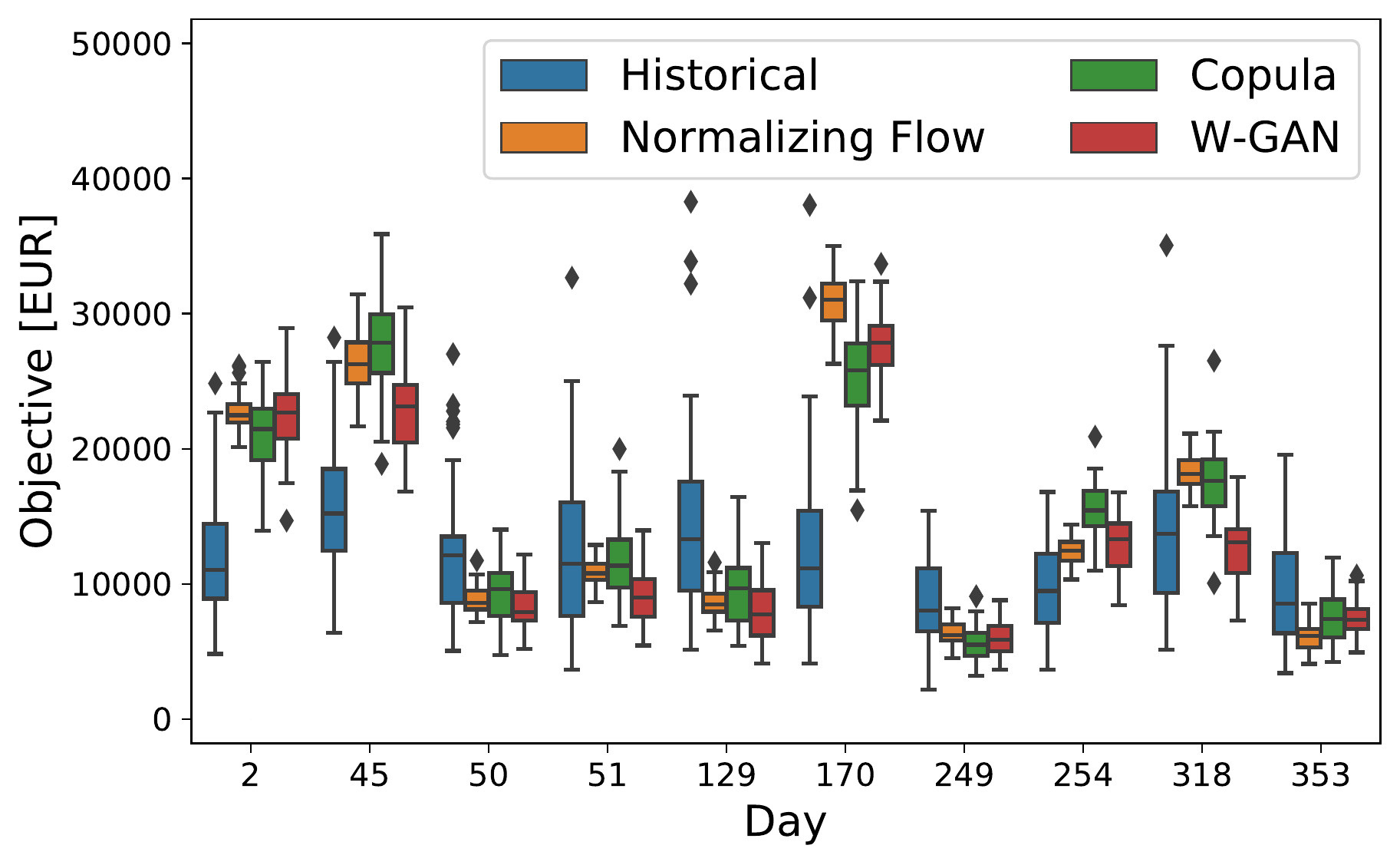}
    \caption{
        Boxplot of the optimal objective function distributions for 50 instances of Problem~\eqref{Prob:Cond_WindProducerProblem} with five scenarios from the historical data set (``Historical''), normalizing flow (``Normalizing Flow''), Gaussian copula (``Copula''), and W-GAN (``W-GAN''), respectively. 
        Boxes indicate quartiles and diamonds indicate outliers \citep{waskom2021seaborn}. 
    }
    \label{fig:Cond_ConsistiencyBoxplot}
\end{figure}
Figure~\ref{fig:Cond_ConsistiencyBoxplot} shows box plots of the optimal objectives using small sets of five scenarios from the historical data set, normalizing flow, Gaussian copula, and W-GAN for ten randomly selected days from 2019.
As an indicator of stability, we look at the spread of optimal objectives, i.e., the height of the box plots, that result from the 50 different scenario sets.
For the ten randomly selected days shown in Figure~\ref{fig:Cond_ConsistiencyBoxplot}, the normalizing flow scenarios result in the lowest spreads of optimal objectives. 
The optimal objectives of the Gaussian copula scenarios show significantly larger spreads than both the normalizing flow and the W-GAN scenarios.
Meanwhile, the randomly selected historical scenario sets lead to by far the largest spreads. 

\begin{table}
    \centering
    \caption{
        Average standard deviation (StD) and average max-min spread (Spread) of the optimal objectives of 50 different instances of Problem~\eqref{Prob:Cond_WindProducerProblem} with 3, 5, 10, 20, and 50 scenarios each generated from normalizing flow, Gaussian copula, and W-GAN over all days in 2019. Best results are marked in \textbf{bold} font.
        }
    \begin{tabularx}{\columnwidth}{Xrrrrr}
    \hline
                            & \# Scenarios     & Historical & Normalizing Flow & Copula & W-GAN \\
                 \hline
StD {[}EUR{]}               & 3                & 6735       & 1726             & 3534   & 2441  \\
                            & 5                & 5025       & 1317             & 2585   & 1883  \\
                            & 10               & 3291       & 930              & 1763   & 1333  \\
                            & 20               & 2253       & 659              & 1244   & 932   \\
                            & 50               & 1404       & \textbf{417}     & 787    & 593   \\ 
\hline
Spread {[}EUR{]}            & 3                & 30527      & 7874             & 17467  & 10719 \\
                            & 5                & 22899      & 5979             & 12285  & 8367  \\
                            & 10               & 14963      & 4214             & 8142   & 5966  \\
                            & 20               & 10138      & 3017             & 5613   & 4184  \\
                            & 50               & 6371       & \textbf{1891}    & 3531   & 2675  \\
\hline
    \end{tabularx}
    \label{tab:Cond_Stability}
\end{table}
Table~\ref{tab:Cond_Stability} shows statistics for the different numbers of scenarios derived over all days in 2019, namely, the average standard deviation and the max-min spread, i.e.,
the difference between the maximum and minimum objective value. 
The results in Table~\ref{tab:Cond_Stability} confirm the observation from Figure~\ref{fig:Cond_ConsistiencyBoxplot} that normalizing flows yield scenarios with the most stable optimal objectives indicated by the lowest standard deviation and the lowest spread.
Furthermore, both standard deviations and spreads decrease for increasing numbers of scenarios for all scenario generation methods. 
Notably, the normalizing flow achieves greater stability with fewer scenarios compared to the other methods. For instance, the standard deviation achieved with just three normalizing flow scenarios is lower than the standard deviations of ten Gaussian copula scenarios and lower than five W-GAN scenarios. 

It appears that the outlier scenarios of the Gaussian copula observed in Figure~\ref{fig:Cond_ExampleScenario} are weighted more heavily in the case of only a few scenarios and, thus, lead to larger differences in the optimal objectives. 
As the normalizing flow shows no extreme outlier scenarios and identifies the overall trends well, there is very little variance in the optimal objective values even for small scenario sets. 
The historical data results in the largest spreads as the data is sampled from the entire spectrum of possible realizations instead of the narrower distributions described by the day-ahead scenario generation models.

\subsection{Obtained profits}
Solving Problem~\eqref{Prob:Cond_WindProducerProblem} maximizes the \textit{expected} profits and yields first-stage decisions that describe a fixed schedule of electricity delivery commitments for the day-ahead market $P^{D}_t$. 
We can compute the \textit{actual} profits by fixing the first stage decisions obtained from the stochastic program and re-optimizing the second stage with the realized electricity production.
To gain statistically relevant results, we compute the \textit{actual} profits obtained in Problem~\eqref{Prob:Cond_WindProducerProblem} for each day in 2019, each time using 100 historical or generated scenarios. 
For comparison, we also compute the profits obtained via the perfect foresight problem, i.e., an ideal single scenario instance of Problem~\eqref{Prob:Cond_WindProducerProblem} based on the actual realization. 

\begin{figure}
    \centering
    \includegraphics[width=0.7\columnwidth]{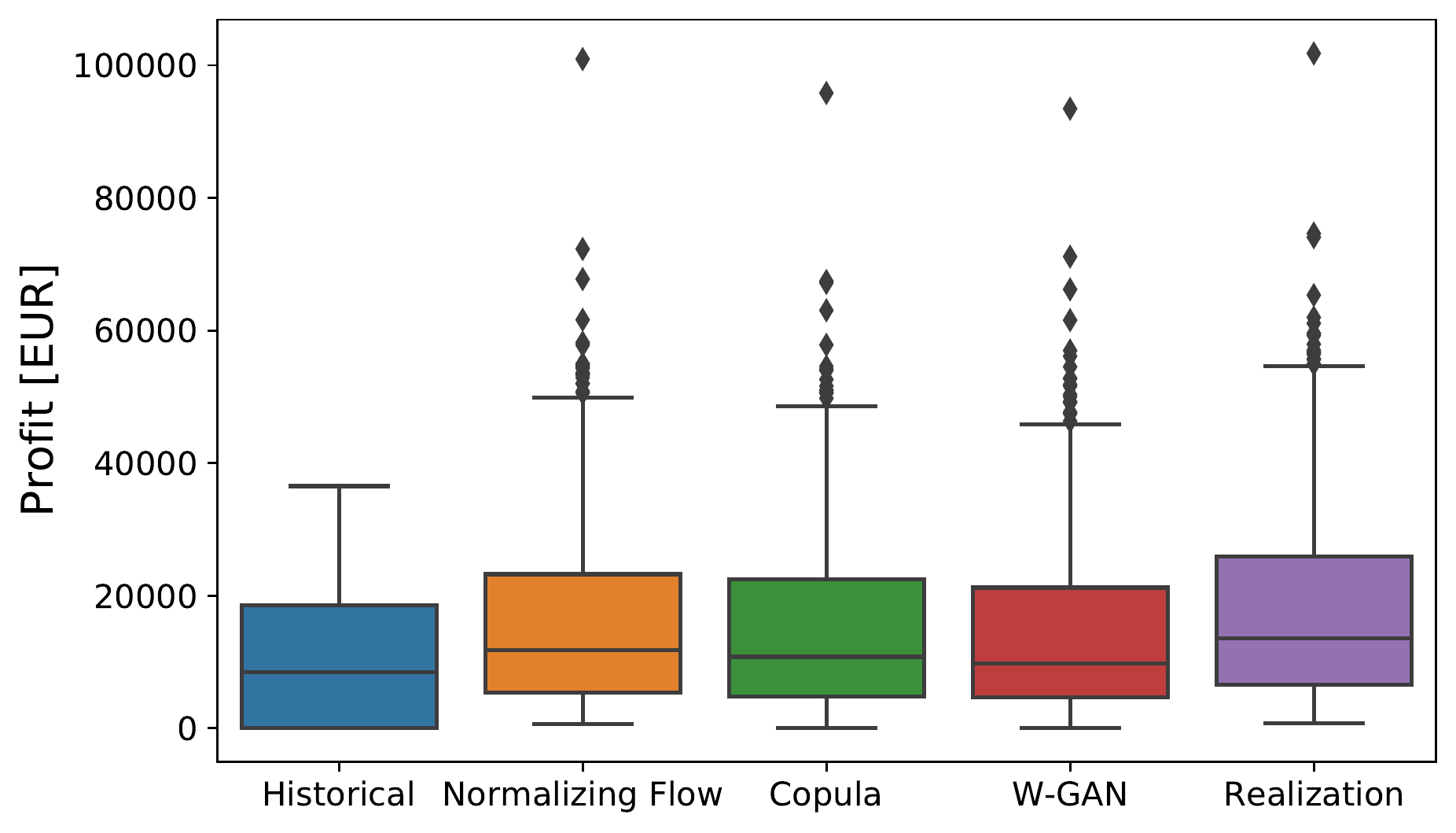}
    \caption{
        Boxplot of profits obtained in 2019 in Problem~\eqref{Prob:Cond_WindProducerProblem}. Each problem uses 100 scenarios generated from the historical data (``Historical''), the normalizing flow (``Normalizing Flow''), the Gaussian copula (``Copula''), and the W-GAN (``W-GAN'') or the realization for the perfect foresight (``Realization''), respectively. 
        Boxes indicate quartiles and diamonds indicate outliers \citep{waskom2021seaborn}.
    }
    \label{fig:Cond_ProfitBoxplot}
\end{figure}
Figure~\ref{fig:Cond_ProfitBoxplot} shows box plots of the distribution of profits obtained by using scenarios from the historical data and the three different generation methods in comparison to the profits obtained from the perfect foresight problem (``Realization''). 
The profits obtained by using the normalizing flow scenarios are the highest, while the Gaussian copula scenarios yield profits between the normalizing flow and the W-GAN.
The profits obtained by using the unconditional historical scenarios are significantly lower than those of all models generating day-ahead scenarios.
Furthermore, the historical scenarios are unconditional and fail to identify the days with high production that appear as high-profit outliers in the case of the other scenario generation methods and the perfect foresight case. 

\begin{table}
    \centering
    \caption{
        Average, average percentage, and maximum perfect information profit gap (PIPG) of \textit{actual} profits obtained over all days in 2019 with 100 scenarios each generated from normalizing flow, Gaussian copula, and W-GAN, respectively. Best results are marked in \textbf{bold} font.}
    \begin{tabularx}{\columnwidth}{Xrrrr}
    \hline
                        & Historical & Normalizing Flow & Copula  & W-GAN   \\ \hline
    Average PIPG [EUR]  & -8381      & \textbf{-1832}   & -2658   & -3485   \\
    Average PIPG [\%]   & -81.5\%    & \textbf{-10.9\%} & -16.6\% & -23.0\%  \\
    Max. PIPG [EUR]     & -70234     & \textbf{-8311}   & -11017  & -17432   \\ 
    \hline
    \end{tabularx}\label{tab:Cond_PIPG_Profits}
\end{table}
To analyze whether the scenario generation methods take advantage of the profit potentials, we define the perfect information profit gap (PIPG) as the difference between the \textit{actual} profits obtained from the stochastic program and the perfect information profit.
Table~\ref{tab:Cond_PIPG_Profits} lists the average and maximum PIPGs in 2019. 
The average PIPGs show that the normalizing flow scenarios yield profits that are on average 6\% and 12\% points closer to the perfect foresight profits compared to the Gaussian copula and the W-GAN, respectively. 
Meanwhile, the historical scenario profits are over 80\% lower on average than the perfect foresight profit. 
The maximum PIPG, i.e., the worst performing days, also show that the normalizing flow scenarios give significantly more profitable results compared to the other generation methods.

\begin{table}
    \centering
    \caption{Average \textit{expected} profits (``Average objective {[}EUR{]}'') of the stochastic program and average \textit{actual} profits (``Average profit {[}EUR{]}'') of 50 different instances of Problem~\eqref{Prob:Cond_WindProducerProblem} with 3, 5, 10, 20, and 50 scenarios each from the historical data and generated from the normalizing flow, the Gaussian copula, and the W-GAN over all days in 2019. 
        }
    \begin{tabularx}{1\columnwidth}{Xrrrrr}
    \hline
                            & \# Scenarios     & Historical & Normalizing Flow & Copula & W-GAN   \\
                 \hline
Average objective {[}EUR{]} & 3                & 12268      & 17336            & 16435  & 16011 \\
Average profit {[}EUR{]}    &                  & 5753       & 16156            & 13527  & 14694 \\
\hline
Average objective {[}EUR{]} & 5                & 11890      & 17195            & 16196  & 15880 \\
Average profit {[}EUR{]}    &                  & 7544       & 16438            & 15197  & 14813 \\
\hline
Average objective {[}EUR{]} & 10               & 10821      & 16974            & 15833  & 15646 \\
Average profit {[}EUR{]}    &                  & 8390       & 16628            & 15738  & 14988 \\
\hline
Average objective {[}EUR{]} & 20               & 10305      & 16893            & 15675  & 15530 \\
Average profit {[}EUR{]}    &                  & 8864       & 16725            & 15880  & 15103 \\
\hline
Average objective {[}EUR{]} & 50               & 10004      & 16829            & 15597  & 15464 \\ 
Average profit {[}EUR{]}    &                  & 9132       & 16792            & 15984  & 15170 \\
\hline
    \end{tabularx}
    \label{tab:Cond_Expected_Vs_Actual_Profits}
\end{table}
Table~\ref{tab:Cond_Expected_Vs_Actual_Profits} lists the average \textit{expected} profits (``Average objective {[}EUR{]}'') and the average \textit{actual} profits (``Average profit {[}EUR{]}'') obtained using scenario sets of size 3, 5, 10, 20, and 50. 
For all scenario generation methods, smaller scenario sets tend to overestimate the \textit{expected} profits. For increasing numbers of scenarios, the \textit{expected} profits and the \textit{actual} profits converge. 
Notably, the normalizing flow scenarios result in the smallest difference between the \textit{expected} and the \textit{actual} profits, particularly for smaller scenario sets of three or five scenarios. 
Furthermore, the average \textit{actual} profits obtained using the normalizing flow are consistently higher than the \textit{actual} profits obtained from any of the other methods. In fact, using just three normalizing flow scenarios results in higher profits than any investigated number of scenarios from any other considered scenario generation method.

The higher profits obtained using the normalizing flow scenarios reflect the findings of Section~\ref{sec:Cond_CaseStudyData}, i.e., the normalizing flow identifies the correct trends and also reflects a diverse distribution.
Meanwhile, the Gaussian copula shows outliers and does not match the distribution, and the W-GAN even struggles to identify the daily trends.
The unconditional historical scenarios do not describe the daily trends and, thus, result in significantly lower profits in the day-ahead scheduling optimization. 
In conclusion, the results shown in Figure~\ref{fig:Cond_ProfitBoxplot}, Table~\ref{tab:Cond_PIPG_Profits}, and Table~\ref{tab:Cond_Expected_Vs_Actual_Profits} highlight that the normalizing flow generates the best scenarios and yields the most profitable bids.   \section{Conclusion}\label{sec:Cond_Conclusion}
The present work considers the scenario generation problem for a day-ahead bidding problem of a wind farm operator to participate in the EPEX spot markets. 
We utilize a data-driven multivariate scenario generation scheme based on conditional normalizing flows to model the distribution of wind capacity factor trajectories with wind speed predictions as conditional inputs. The generated scenarios are specifically tailored to stochastic optimization problems concerning the time frame between 00:00\,am and 11:59\,pm.

We analyze the normalizing flow scenarios in comparison to randomly selected historical data and scenarios generated from other more established methods, namely, Gaussian copulas and Wasserstein generative adversarial networks (W-GANs), and compare them to the actually realized power generation in 2019.
The historical scenarios reflect the overall distribution of realizations well but fail to identify daily trends and show large variations independent of the investigated day.
Among the conditional scenario generation methods, the normalizing flow scenarios best reflect the realized power generation trends and their distributions while also displaying a diverse set of possible realizations. 
Meanwhile, the Gaussian copula results in uncharacteristic outliers, and the W-GAN struggles to identify the main trends of the realizations. Furthermore, both Gaussian copula and W-GAN result in skewed distributions. 

To assess their value for stochastic programming, the scenarios are applied in a stochastic programming case study that aims to set bids for electricity sales on the day-ahead market. 
First, we investigate the stability defined by \cite{Kaut2003Evaluation}. Here, the normalizing flow yields scenarios that result in the most stable optimal objective values, even for small scenario sets. 
In particular, the normalizing flow scenarios result in the smallest standard deviation and the smallest spread in the optimal objective values.
The analysis of the actual profits obtained on all days of 2019 shows a significant advantage of using day-ahead scenarios that are specifically tailored to the investigated day. Using randomly selected historical scenarios results in an average perfect information profit gap (PIPG) of over 80\%, while the conditional scenario generation methods return PIPGs between 10 to 23\%. 
The bids placed using the normalizing flow scenarios obtain the highest profits and have the lowest PIPG, i.e., the profits are closest to the perfect foresight profits. Notably, even small scenario sets of only three normalizing flow scenarios result in higher actual profits than any investigated number of scenarios from the other considered scenario generation methods.

In conclusion, utilizing conditional, day-specific scenarios in day-ahead scheduling problems leads to significantly more profitable decisions compared to relying on unconditional historical data. 
Furthermore, the conditional normalizing flow model yields high-quality scenarios that result in highly profitable solutions for stochastic programs and stable results, even for small scenario sets.
Therefore, we argue that normalizing flow scenarios have a high potential for scheduling problems that cannot be solved with a large scenario set. 
 
\section*{Acknowledgements}
\noindent We would like to thank Marcus Vo{\ss} (Technical University of Berlin, Distributed Artificial Intelligence Laboratory) for his valuable input on Copula methods, scenario evaluation, and the supervision of L. Paeleke.
This work was performed as part of the Helmholtz School for Data Science in Life, Earth and Energy (HDS-LEE) and received funding from the Helmholtz Association of German Research Centres.

  \appendix
  \bibliographystyle{apalike}
\renewcommand{\refname}{Bibliography}
  \bibliography{extracted.bib}

\end{document}


\thispagestyle{firststyle}

  \begin{center}
    \begin{large}
      \textbf{\mytitle}
    \end{large} \\
    \myauthor
  \end{center}

  \vspace{0.5cm}

  \begin{footnotesize}
    \affil
  \end{footnotesize}

\section{Daily quantiles}
\begin{figure*}[h!]
    \centering
    \includegraphics[width=\textwidth]{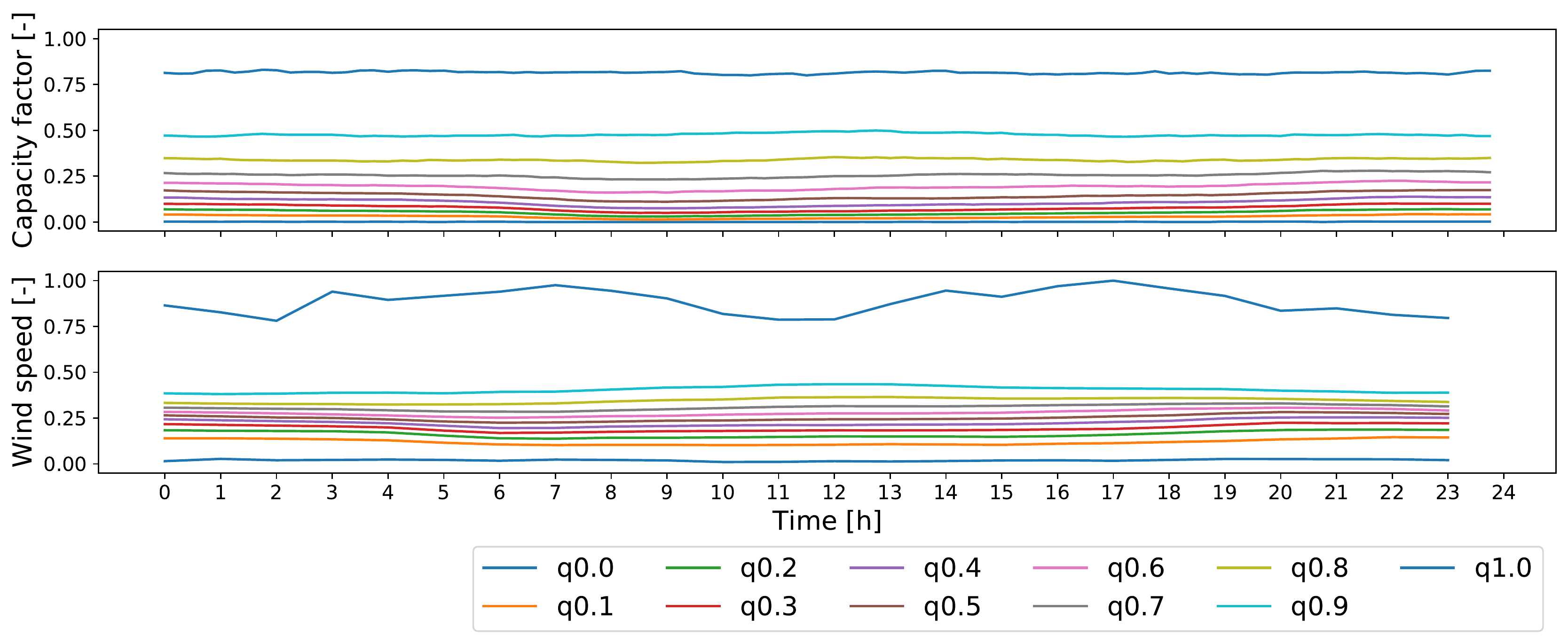}
    \caption{
        Quantile trajectories of wind capacity factor (top) and scaled wind speed (bottom). 
        Quantiles are computed for every time step individually. 
        q0.0 and q1.0 refer to the lowest and the highest values in the data set, respectively. 
        }
    \label{fig:Cond_Supp_QuantileTrajectories}
\end{figure*}

Figure~\ref{fig:Cond_Supp_QuantileTrajectories} plots the quantiles of the capacity factor and wind speed distributions over the scenario length of one day. Neither Capacity factor nor wind speed show any significant daily trends. The distributions are very broad, motivating the use of conditional data.

\section{Scenario generation model implementation}
Both normalizing flow and W-GAN use both Gaussian samples and conditional inputs as direct inputs to their respective ANN structures \citep{chen2018model, zhang2019scenario}. 
The Gaussian copula uses linear quantile-regression based on the conditional inputs to estimate the inverse cumulative distribution function (CDF), which is then used to transform the Gaussian samples \citep{pinson2009probabilistic}.
For more information on scenario generation using Gaussian copulas and W-GANs, we refer to \cite{pinson2009probabilistic} and \cite{chen2018model}, respectively.

Due to the significant dimensionality reduction from 96 to 18 dimensions, a small RealNVP \citep{dinh2016realNVP} model is sufficient. The employed model uses four affine coupling layers with fully connected conditioner models with 2 hidden layers with 9 neurons each. 
The Gaussian copula was implemented using the linear quantile regression in \cite{seabold2010statsmodels} and the required inverse CDF is estimated using linear interpolation with 20 intervals. 
The model structures of the generator and critic used for the W-GAN are shown in Table~\ref{tab:Cond_WGAN_Layers}. 
Both the normalizing flow model and the W-GAN are implemented in TensorFlow 2.5.0~\citep{tensorflow2015}. The PCA is computed using the scikit-learn library \citep{scikitlearn}.

\begin{table}
    \centering
    \caption{
        Layers of generator and critic for W-GAN \citep{arjovsky2017wasserstein}. Attributes of layer types are: 
            Linear (fully connected): Number of nodes, 
            Conv (1D convolutional layer): (Number of filters, filter size, strides, padding), 
            Conv-T (1D convolutional layer transpose): (Number of filters, filter size, strides, padding), and 
            Reshape: (output dim 1, output dim 2, \dots). 
            The W-GAN is trained for 2000 epochs and the sampling dimensionality is 20.
            }

    \begin{tabular}{lrl}
\hline
        Generator layers   & Attributes    & Activation \\
        \hline
        Linear  & 96            & ReLU          \\ 
        Linear  & $96\cdot 12$  & ReLU          \\
        Reshape & (96,12)       & -             \\
        Conv-T  & (12,3,1,1)    & ReLU          \\
        Conv-T  & (1,3,1,1)     & tanh          \\\hline
        \\ 
\hline
        Critic layers   & Attributes    & Activation \\ 
        \hline
        Conv    & (12,3,1,1)    & LeakyReLU     \\
        Conv    & (4,3,1,1)     & LeakyReLU     \\
        Flatten & -             & - \\
        Linear  & 96            & LeakyReLU \\
        Linear  & 1             & - \\ 
        \hline
    \end{tabular}
    \label{tab:Cond_WGAN_Layers}
\end{table}     
  
  \appendix
  \bibliographystyle{apalike}
\renewcommand{\refname}{Bibliography}
  \bibliography{extracted.bib}